\documentclass[conference,a4paper]{IEEEtran}
\usepackage{xcolor}
\usepackage{balance}

\ifCLASSINFOpdf
  \usepackage[pdftex]{graphicx}
\else
  \usepackage[dvips]{graphicx}
\fi
\usepackage{amsmath}
\usepackage{bm,amssymb}
\ifCLASSOPTIONcompsoc
 \usepackage[caption=false,font=normalsize,labelfont=sf,textfont=sf]{subfig}
\else
 \usepackage[caption=false,font=footnotesize]{subfig}
\fi

\usepackage{stfloats}
\usepackage{cite}
\usepackage{calc,ifthen,intcalc}
\usepackage{tikz}
	\usetikzlibrary{arrows,shapes,positioning,calc,fadings,patterns,decorations}
	\usetikzlibrary{decorations.pathreplacing,decorations.pathmorphing}			
	\usetikzlibrary{arrows,arrows.meta}
	\usetikzlibrary{shapes.geometric}

\usepackage{booktabs}

\newcommand{\bs}[1]{\bm{\mathrm{#1}}}
\newcommand{\abs}[1]{\left|#1\vphantom{f}\right|}
\renewcommand\Re{\operatorname{Re}}
\renewcommand\Im{\operatorname{Im}}
\newcommand\norm[1]{\ensuremath{\left\|#1\right\|}}
\renewcommand{\d}{\ensuremath{\partial}}

\newcommand{\tr}{\ensuremath{\operatorname{tr}}}

\newcommand\figref[1]{Fig.~\ref{#1}}

\usepackage{lipsum}

\begin{document}
\bstctlcite{IEEEexample:BSTcontrol}
\title{Investigating Sparse Reconfigurable Intelligent Surfaces (SRIS) via Maximum Power Transfer Efficiency Method Based on Convex Relaxation}


\author{\IEEEauthorblockN{
Hans-Dieter Lang\IEEEauthorrefmark{1},
Michel A. Nyffenegger\IEEEauthorrefmark{1},
Heinz Mathis\IEEEauthorrefmark{1}, and
Xingqi Zhang\IEEEauthorrefmark{2}}
\IEEEauthorblockA{\IEEEauthorrefmark{1}OST -- Eastern Switzerland University of Applied Sciences, Rapperswil, SG, Switzerland. Email: hansdieter.lang@ost.ch}
\IEEEauthorblockA{\IEEEauthorrefmark{2}School of Electrical and Electronic Engineering, University College Dublin, Dublin, Ireland.}}



\maketitle

\begin{abstract}

Reconfigurable intelligent surfaces (RISs) are widely considered to become an integral part of future wireless communication systems. Various methodologies exist to design such surfaces; however, most consider or require a very large number of tunable components. This not only raises system complexity, but also significantly increases power consumption. Sparse RISs (SRISs) consider using a smaller or even minimal number of tunable components to improve overall efficiency while maintaining sufficient RIS capability. The versatile semidefinite relaxation-based optimization method previously applied to transmit array antennas is adapted and applied accordingly, to evaluate the potential of different SRIS configurations. Because the relaxation is tight in all cases, the maximum possible performance is found reliably. Hence, with this approach, the trade-off between performance and sparseness of SRIS can be analyzed. Preliminary results show that even a much smaller number of reconfigurable elements, e.g. only 50\%, can still have a significant impact.

%This allows analyzing the trade-off between performance and sparsity of SRISs, with preliminary results showing that even a considerably smaller number of reconfigurable elements, e.g., just 50\%, can still have a substantial impact.
\end{abstract}

\vskip0.5\baselineskip

% \fi
% \IEEEpeerreviewmaketitle

\vspace*{4pt}
\section{Introduction}
Reconfigurable intelligent surfaces (RISs) continue to gain interest in many research areas including applied electromagnetics, antennas \& propagation, and mobile communications~\cite{hu18,huang19,renzo20,wu20}. They are widely considered to play a key role in future wireless communication technologies and networks, such as "beyond 5G", 6G, and onwards \cite{huang20,saad20,liu21}. 

The main idea builds on reconfigurable reflectarrays~\cite{hum14}: RIS divert otherwise unused ("spilled") wireless signals back into the system, to enhance coverage and throughput of wireless links. 
In order for the reflection to be useful, these reflectors must be both reconfigurable and "intelligent", i.e., they must be controlled to optimize the overall performance for each link, taking into account all the users involved as well as the entire channel, including all other scatterers. 
Such smart reflectors should reduce the need for additional base stations and save both costs and energy. 
Therefore, the use of RIS should not only improve system performance per cost, but also overall energy efficiency and sustainability.

The last point is one of many challenges associated with this technology.
%One of the many challenges faced with this technology is power efficiency. 
The premise is that such RISs constitute (much) more power-efficient network components than additional fully functional base stations. 
RIS are considered to be low-power and low-cost \cite{taghvaee22}. However, they are not entirely passive, as they must contain tunable elements. Various tuning mechanisms have been considered \cite{tsilipakos20}. 
For applications below optical and mmWave frequency ranges, most importantly biased varactor diodes are used as tunable capacitors. And while each varactor diode can be considered low-power, a large number of them  can still add up to a considerable amount of power required for the operation of the RIS (in addition to the power consumption of the central control engine).%, so that the reflection (or transmission) can be steered onto the desired target.

\vspace*{4pt}
\section{Sparse RIS (SRIS)}% for increased energy efficiency}
%\textcolor{red}{Design procedures include...} holography \cite{huang20}
%Other problems include network optimization \cite{renzo22}

A large number of tunable elements clearly increases the number of degrees of freedom for accuracy and focusing of the reflected beam (possibly to many users at the same time, at different frequencies or in different time slots). However, this not only increases the system complexity when it comes to setting all these elements correctly (i.e., with an optimization algorithm), but also the total power consumption of the RIS. Therefore, one of the criticisms is that such RIS are not necessarily simpler and at least not much more energy efficient than additional base stations.

Hence, the goal is to reduce the number of tunable elements required, to improve the energy efficiency of the RIS.
The result is an RIS with sparsely populated tunable elements, a so-called Sparse RIS (SRIS). 
Sparse RISs (SRISs) consider using a smaller or even minimal number of such tunable components to improve efficiency while maintaining sufficient RIS capability. 
At this stage, it is unclear which strategy will best accomplish this. 
The most commonly used techniques to design RISs, see, e.g., \cite{rahmat19,huang20},  seem not directly applicable. 

This is where the optimization procedure adopted from array antennas, as described in the following, comes in: Both to assess the performance of a particular RIS design as well as to investigate the effects of reducing the number of tunable elements on an RIS, e.g., by clustering or simply removing some tunable elements, a method is desired that give the maximum enhancement for a given use case. 

Note that, beyond the intrinsic sparsity of RISs, employing sparsity also in the deployment of these SRISs further improves sustainability, by minimizing both manufacturing costs and resource consumption, as well as overall energy consumption. However, this is not a topic for this work.

%\subsection{Optimization for SRIS Design \& Assessment}

%\vspace*{4pt}
\section{Optimization Procedure}

The following optimization procedure is based on a framework originally developed for wireless power transfer (WPT) in the near-field \cite{lang17,lang22}, which has subsequently been adapted for far-field applications \cite{lang17ant}. More recently also sidelobe and backlobe constraints have been added  \cite{nyffenegger22}, which are beyond the scope of this paper, but are generally interesting for (S-)RIS analysis as well.

The overall objective of the method is to find the optimal feasible currents in all antenna elements that minimize the total transmit power required to obtain unity power in the receive antenna. %Feasibility pertains to enforcing KVL as well as passivity of all RIS elements. 
For this purpose, the entire SRIS setup is simulated as a multiport network to obtain its impedance matrix. The required Kirchhoffs Voltage Law (KVL) and power equations are then all formulated using standard circuit network analysis based on that matrix. For the SRIS to be passive, all average transferred powers at RIS element ports are required to be zero, while the transmit power is minimized.
%Using standard circuit network analysis, the corresponding equations are formulated and solved to maximize the received power.

%Additional power constraints allow impedance matching requirements to be taken into account as well \cite{lang18eucap}. However, they are not needed here, since the reflectarray is usually sufficiently far away from the feed antenna, such that detuning can be neglected.
%In general, this optimization problem is nonconvex; colloquially often referred to as ``unsolvable'', for a large number of variables (i.e., the number of RIS elements). However, a reformulation %in terms of a generalized current matrix, in combination with a typical semidefinite relaxation, transforms the problem into a convex semidefinite program (SDP) \cite{boyd04}. Such programs are reliably and efficiently solvable, even for a (reasonably) large number of variables, using dedicated algorithms, such as SDPT3 \cite{tutuncu03}. % in Matlab.

Note that this method is by no means limited to surfaces considering reflection only; it may just as well be applied to analyze transmitting surfaces or surfaces that simultaneously transmit and reflect, e.g., STAR-RIS \cite{xu21,liu21star}. Furthermore, the method can also be used to assess the performance of a particular RIS design, when using reactively loaded RIS elements (instead of geometrically adjusting the elements to achieve essentially the same reactive behavior).

\subsection{Prerequisites}% in Impedance Matrix Form}

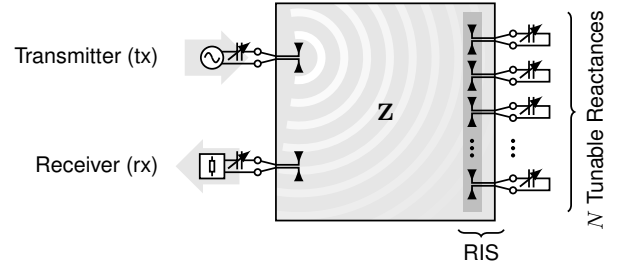
\begin{figure}[!ht]\centering
		%\vspace*{-12mm}
	\begin{tikzpicture}[>=latex,scale=1.2,line width=0.66pt]\sffamily\small
			
			\fill[gray!20!white] (-2.6,0.6)+(-0.2,-0.2)-|+(0.2,-0.2-0.1)--+(0.5,0)--+(0.2,0.2+0.1)|-+(-0.2,0.2)--cycle;
			\fill[gray!20!white] (-2.4,-0.6)+(0.2,-0.2)-|+(-0.2,-0.2-0.1)--+(-0.5,0)--+(-0.2,0.2+0.1)|-+(0.2,0.2)--cycle;
		
			%\path (-1.5,-1.3)--(3,-1.3);
			%\path (-1.5,-1.3)rectangle(3,1.3);
			\fill[gray!20!white] (0,0)+(0.6-2.4,-1.2)rectangle+(0.6,1.2);
			
			\begin{scope}
			\clip  (0,0)+(0.6-2.4,-1.2)rectangle+(0.6,1.2);
			%\draw[decorate,decoration={expanding waves,angle=80},line width=2pt,white,path fading=east] (-1.55,0.6)--+(0:3);
			\def\xtx{-1.55}
			\draw[white,line width=3pt] (\xtx,0.6)+(100:0.25)arc(100:-100:0.25);
			\draw[white,line width=3pt,opacity=0.8] (\xtx,0.6)+(100:0.5)arc(100:-100:0.5);
			\draw[white,line width=3pt,opacity=0.7] (\xtx,0.6)+(100:0.75)arc(100:-100:0.75);
			\draw[white,line width=3pt,opacity=0.6] (\xtx,0.6)+(100:1)arc(100:-100:1);
			\draw[white,line width=3pt,opacity=0.5] (\xtx,0.6)+(100:1.25)arc(100:-100:1.25);
			\draw[white,line width=3pt,opacity=0.4] (\xtx,0.6)+(100:1.5)arc(100:-100:1.5);
			\draw[white,line width=3pt,opacity=0.3] (\xtx,0.6)+(100:1.75)arc(100:-100:1.75);
			\draw[white,line width=3pt,opacity=0.2] (\xtx,0.6)+(100:2)arc(100:-100:2);
			\draw[white,line width=3pt,opacity=0.1] (\xtx,0.6)+(100:2.25)arc(100:-100:2.25);
			\end{scope}
			
			\draw(0,0)+(0.6-2.4,-1.2)rectangle+(0.6,1.2);
			
			%\fill[pattern=north west lines,line width=0.25pt,pattern color=gray] (0.45,-1.1)rectangle(0.25,1.1);
			\fill[gray!50!white] (0.45,-1.1)rectangle(0.25,1.1);
			\def\dy{0.066}

			\foreach \ysource in {-0.8,0,0.4,0.8}{
				\draw (0.8,\ysource-\dy)--(1.2,\ysource-\dy)--(1.2,\ysource+\dy)--(1+0.025,\ysource+\dy)
				(+1+0.025,\ysource+\dy-0.1)--(1+0.025,\ysource+\dy+0.1)
				(+1-0.025,\ysource+\dy-0.1)--(1-0.025,\ysource+\dy+0.1)
				(+1-0.025,\ysource+\dy)--(0.8,\ysource+\dy);
				
				\draw (0.8,\ysource+\dy)--+(-0.2,-0.05);
				\draw (0.8,\ysource-\dy)--+(-0.2,0.05);
				%\draw[path fading=west] (0.6,\ysource+\dy)--+(-0.25,0);
				%\draw[path fading=west] (0.6,\ysource-\dy)--+(-0.25,0);
				\draw[-<] (0.6,\ysource+\dy-0.05)-|+(-0.25,0.15);
				\draw[-<] (0.6,\ysource-\dy+0.05)-|+(-0.25,-0.15);
				\draw[fill=white] (0.8,\ysource+\dy)circle(1pt);
				\draw[fill=white] (0.8,\ysource-\dy)circle(1pt);
				\draw[line width=0.7pt,->] (1,\ysource+0.075)+(-0.125,-0.10)--+(0.125,0.10);
			}
			
			\begin{scope}[xshift=-1.2cm,xscale=-1]
			\foreach \ysource in {-0.6,0.6}{
				\draw (0.8,\ysource-\dy)--(1.2,\ysource-\dy)
				(1.2,\ysource+\dy)--(1+0.025,\ysource+\dy)
				(+1+0.025,\ysource+\dy-0.1)--(1+0.025,\ysource+\dy+0.1)
				(+1-0.025,\ysource+\dy-0.1)--(1-0.025,\ysource+\dy+0.1)
				(+1-0.025,\ysource+\dy)--(0.8,\ysource+\dy);
				
				\draw (0.8,\ysource+\dy)--+(-0.2,-0.05);
				\draw (0.8,\ysource-\dy)--+(-0.2,0.05);
				%\draw[path fading=west] (0.6,\ysource+\dy)--+(-0.25,0);
				%\draw[path fading=west] (0.6,\ysource-\dy)--+(-0.25,0);
				\draw[-<] (0.6,\ysource+\dy-0.05)-|+(-0.25,0.15);
				\draw[-<] (0.6,\ysource-\dy+0.05)-|+(-0.25,-0.15);
				\draw[fill=white] (0.8,\ysource+\dy)circle(1pt);
				\draw[fill=white] (0.8,\ysource-\dy)circle(1pt);
				\draw[line width=0.7pt,->] (1,\ysource+0.075)+(0.125,-0.10)--+(-0.125,0.10);
			}
			\end{scope}
			
			\draw[fill=white] (-2.5,-0.6)+(-0.135,-0.135)rectangle+(0.135,0.135);
			\draw (-2.5,-0.6)+(0,-0.11)--+(0,0.11);
			\draw[fill=white] (-2.5,-0.6)+(-0.025,-0.06)rectangle+(0.025,0.06);
			
			\def\xsource{-2.5}
			\foreach \ysource in {0.6}{
				\draw[fill=white,line width=0.66pt] (\xsource,\ysource)circle(3.5pt);
				\draw[line width=0.66pt] (\xsource,\ysource)..controls(\xsource+0.027,-0.07+\ysource)and(\xsource+0.054,-0.07+\ysource)..+(+0.081,0);
				\draw[line width=0.66pt] (\xsource,\ysource)..controls(\xsource-0.025,0.07+\ysource)and(\xsource-0.054,0.07+\ysource)..+(-0.081,0);
			}
				
			\node[inner sep=2pt] (z0) at(0.6-1.2,0){$\bs{Z}$};

			\draw[decorate,decoration=brace] (0.7,-1.3)--node[below=2.5pt,rotate=0,align=center,scale=0.9]{RIS}(0.2,-1.3);
			\draw[decorate,decoration=brace] (1.4,1.1)--node[below=3pt,rotate=90,yshift=-1.5mm,align=center,scale=0.9]{$N$ Tunable Reactances}(1.4,-1.1);
			\fill (0.6+0.2,-0.4+0.08)circle(0.66pt);
			\fill (0.6+0.2,-0.4)circle(0.66pt);
			\fill (0.6+0.2,-0.4-0.08)circle(0.66pt);
			\fill (0.6-0.25,-0.4+0.08)circle(0.66pt);
			\fill (0.6-0.25,-0.4)circle(0.66pt);
			\fill (0.6-0.25,-0.4-0.08)circle(0.66pt);
			
			\node[anchor=east,scale=0.9] at(-3,0.6){Transmitter (tx)};
			\node[anchor=east,scale=0.9] at(-3,-0.6){Receiver (rx)};
		\end{tikzpicture}
		\vspace*{-1mm}
		\caption{The multiport network to model the wireless link between transmitter, receiver and the RIS, where the goal is to maximize power transfer from the transmitter to the receiver using optimally tuned reactances.}
		\label{fig:problem}
\end{figure}

\begin{figure*}[!htb]\centering
	\begin{tikzpicture}[>=latex]\small\sffamily
		\node[inner sep=0pt] at(0,0){\includegraphics[clip,trim={0 0 0 0},width=0.95\textwidth]{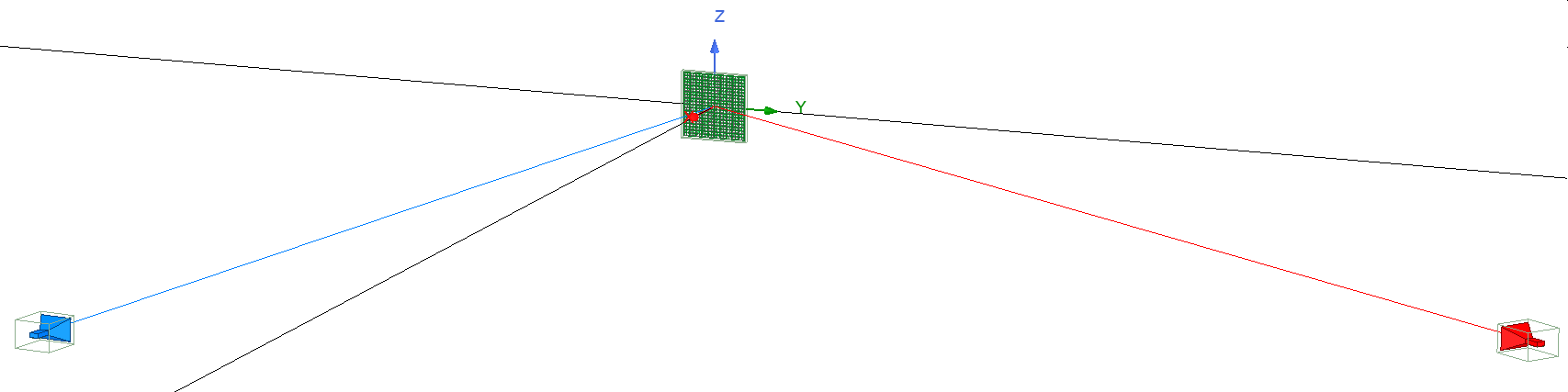}};
		\begin{scope}[xshift=5mm]
		\begin{scope}
		\clip (-1.5,0.5)--(1.5,0.15)--(2.3,-0.15)--(2.3,-3.2)--(-1.5,-2.8)--cycle;
		\node at(0.5,-1.5){\includegraphics[clip,trim={25mm 0 45mm 0},width=4cm]{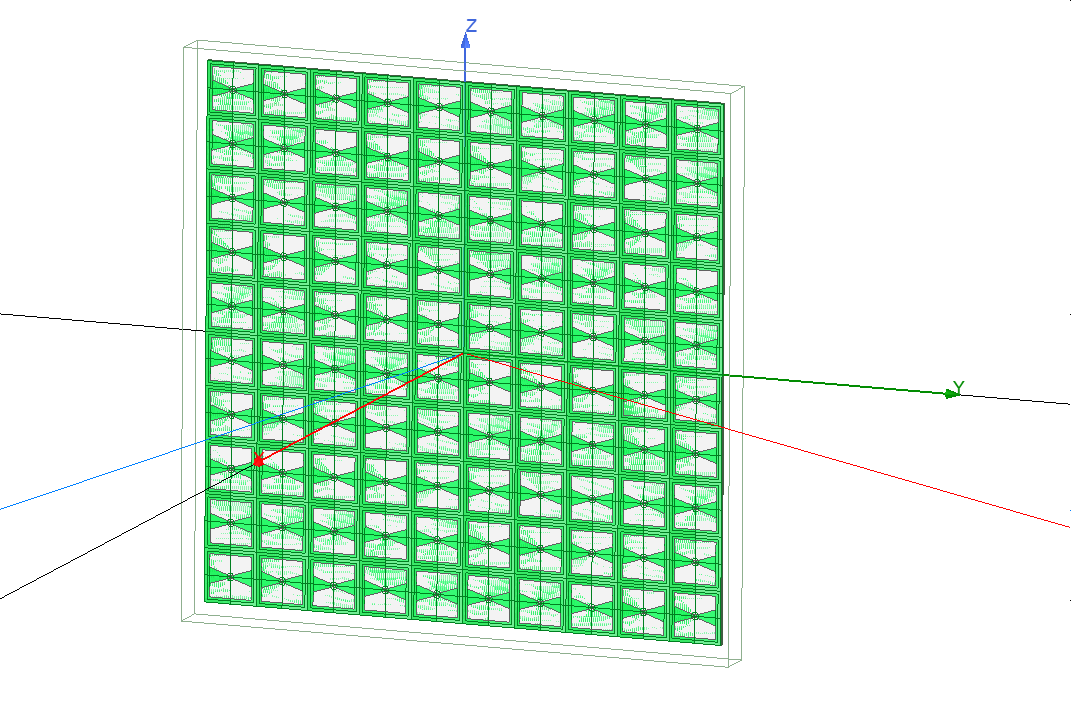}};
		\end{scope}
		\draw[<->] (-0.95,-2.2)--node[above,scale=0.85,rotate=-3]{40\,cm}(1.6,-2.42);
		\draw[line width=0.3pt] (1.63,-2.9)--coordinate[pos=0.8](a)+(-4.6:0.8);
		\draw[line width=0.3pt] (1.63,-0.22)--coordinate[pos=0.8](b)+(-4.6:0.8);
		\draw[<->] (b)--node[above,scale=0.85,rotate=90]{40\,cm}(a);
		\end{scope}
		
		\draw[line width=0.25pt] (0.5,-1.5)+(-0.695,0.4)coordinate(a)--+(-0.695,0.135)coordinate(d)--+(-0.43,0.12)coordinate(c)--+(-0.43,0.385)coordinate(b)--cycle;
		\draw[line width=0.15pt,gray] (a)--(-3.98,-0.61);
		\draw[line width=0.15pt,gray] (b)--(-1.97,-0.74);
		\draw[line width=0.15pt,gray] (c)--(-1.97,-2.85);
		\draw[line width=0.15pt,gray] (d)--(-3.98,-2.85+0.15);
		
		%\clip (-1.5,-1.5)circle(1.3cm);
		\begin{scope}
		\clip (-2.85,-1.75)+(-1.15,1.16)--+(0.88,1.02)--+(0.88,-1.11)--+(-1.15,-0.95)--cycle;
		\node at(-2.85,-1.75){\includegraphics[clip,trim={25mm 0 25mm 0},width=2.8cm]{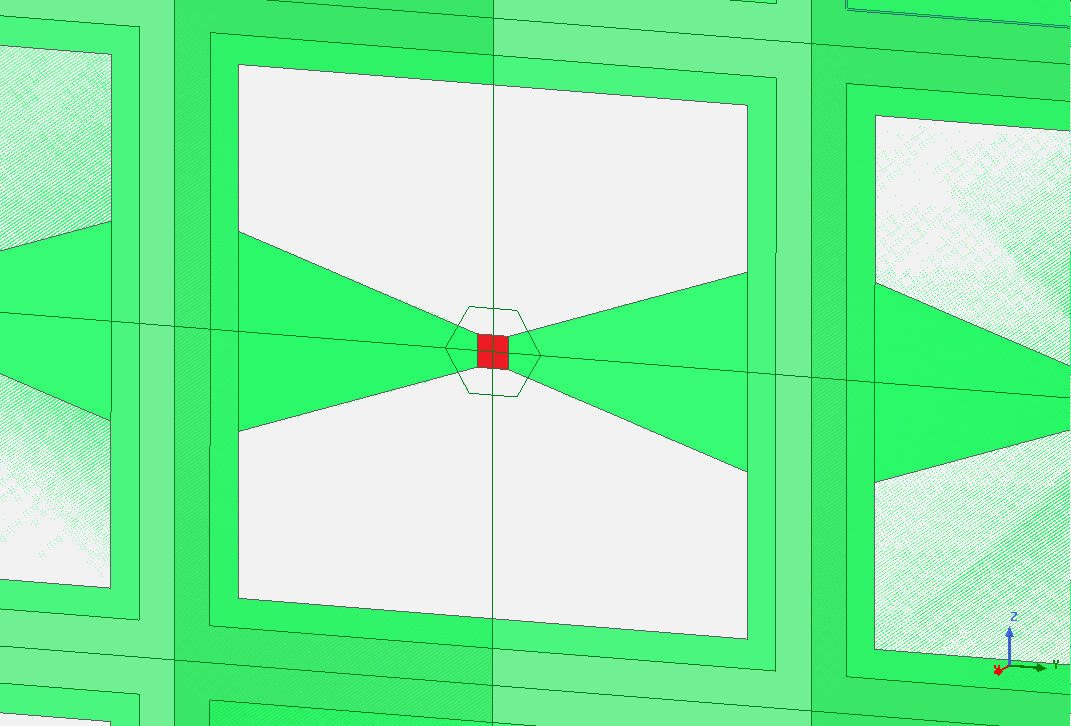}};
		\end{scope}

		\draw[line width=0.3pt] (-2.2,-2.62)--coordinate[pos=0.8](aa)+(-4.6:0.8);
		\draw[line width=0.3pt] (-2.2,-0.93)--coordinate[pos=0.8](bb)+(-4.6:0.8);
		
		\draw[<->] (-3.8,-2.3)--node[above,scale=0.85,rotate=-3]{32\,mm}(-2.19,-2.42);
		\draw[<->] (bb)--node[above,scale=0.85,rotate=90]{32\,mm}(aa);
		\draw (-3,-1.72)--+(-1.1,-0.15)node[left,scale=0.85]{Lumped port};
		\draw (-1.0,1.2)--+(-0.5,0.25)node[left,scale=0.85]{(S-)RIS};
		\node[cyan!66!blue,scale=0.85,align=center] at(-8.1,-2.2){Transmitter\\antenna};
		\node[red,scale=0.85,align=center] at(8.2,-2.25){Receiver\\antenna};
	\end{tikzpicture}
	\vspace*{-4mm}
	\caption{Simulation setup of the test cases in ANSYS HFSS: The transmitter antenna (blue, left) is at an angle 10° off broadside, the receiving antenna (red, right) is at various angles $\alpha$ off broadside and 1° below the $xy$-plane. The RIS consists of $10\times10=100$ vertically polarized printed dipoles.}
	\label{fig:hfss_setup}
	\vspace*{-2mm}
\end{figure*}
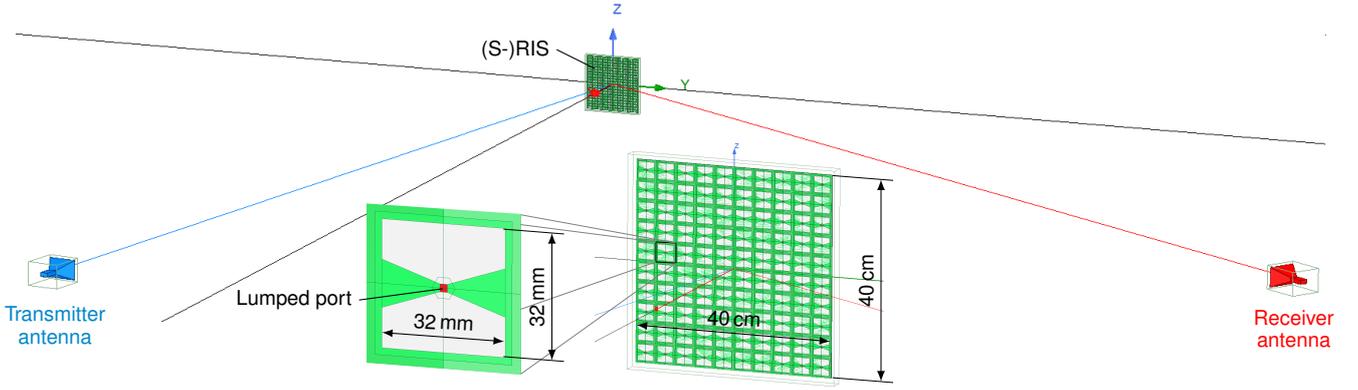

The overview of the problem is visualized in \figref{fig:problem} and an example of such a setup, which is going to serve as test case, is shown in \figref{fig:hfss_setup}. 
In essence, it all boils down to a multiport network problem of transferring maximum power from the transmit port (tx) to the receiver port (rx) with the help of the rest of the network when all other ports are only connected by lossless reactances.

%Since the optimization procedure is derrived in detail in \cite{} and \cite{}, with two small typos corrected in \cite{}, here only the most important are elaborated.

%\subsection{Prerequisites}
	Let the impedance matrix $\bs{Z}$ of the entire (unoptimized) wireless link consisting of the transmit antenna, RIS, and receive antenna be partitioned according to 
	\begin{equation}
		\bs{Z}=
		\begin{bmatrix}
			\tikz{\path[] (-0.2,0.4)--(0.2,0.4);
			\node[fill=gray!35!white,opacity=0.85,minimum width=5mm,minimum height=5mm,inner sep=0pt] at(0,0){$z_t$};}
			\hspace*{-1.25mm}
			&
			\hspace*{-1.25mm}
				\tikz{\node[fill=gray!35!white,opacity=0.85,minimum width=10mm,minimum height=5mm,inner sep=0pt] at(0,0){$\bs{z}_{ts}$};}
			\hspace*{-1.25mm}
			&
			\hspace*{-1.25mm}
			\tikz{\node[fill=gray!35!white,opacity=0.85,minimum width=5mm,minimum height=5mm,inner sep=0pt] at(0,0){${z}_{tr}$};}
			\hspace*{-0.0mm}\\[-0.25mm]
			\tikz{\node[fill=gray!35!white,opacity=0.85,minimum width=5mm,minimum height=10mm,inner sep=0pt] at(0,0){$\bs{z}_{ts}^T$};}
			\hspace*{-1.25mm}
			&
			\hspace*{-1.25mm}
			\tikz{%\draw(-0.4,0.75)--(0.5,0.75);
				\node[fill=gray!35!white,opacity=0.85,minimum width=10mm,minimum height=10mm,inner sep=0pt] at(0,0){$\bs{Z}_s$};}
			\hspace*{-1.25mm}
			&
			\hspace*{-1.25mm}
			\tikz{\node[fill=gray!35!white,opacity=0.85,minimum width=5mm,minimum height=10mm,inner sep=0pt] at(0,0){$\bs{z}_{sr}$};}
			\hspace*{-0.0mm}\\[-0.25mm]
			\tikz{\node[fill=gray!35!white,opacity=0.85,minimum width=5mm,minimum height=5mm,inner sep=0pt] at(0,0){${z}_{tr}$};}
			\hspace*{-1.25mm}
			&
			\hspace*{-1.25mm}
			\tikz{\node[fill=gray!35!white,opacity=0.85,minimum width=10mm,minimum height=5mm,inner sep=0pt] at(0,0){$\bs{z}_{sr}^T$};}
			\hspace*{-1.25mm}
			&
			\hspace*{-1.25mm}
			\tikz{\node[fill=gray!35!white,opacity=0.85,minimum width=5mm,minimum height=5mm,inner sep=0pt] at(0,0){$z_r^{\vphantom{T}}$};}
			\hspace*{-0.0mm}
		\end{bmatrix}
%		\bs{Z}=\begin{bmatrix}
%			\bs{Z}_t & \bs{Z}_{ts} & \bs{z}_{tr}\\[0.2em]
%			\bs{Z}_{ts}^T &  \bs{Z}_s & \bs{z}_{sr}\\[0.2em]
%			\bs{z}_{tr}^{T} & \bs{z}_{sr}^T & z_r\end{bmatrix}
			\in\mathbb{C}^{(N+2)\times (N+2)}\;,
			\label{eq:Z}
	\end{equation}
where $N$ stands for the number of antennas (ports) of the RIS which contain tunable elements. The subscripts $t$, $s$, and $r$ refer to the \underline{t}ransmitter antenna port, the $N$ \underline{s}urface antennas and the \underline{r}eceiver antenna, respectively. The subscript combinations $ts$, $tr$, and $sr$ stand for the coupling of the respective elements. 
	Due to reciprocity, the matrix is symmetric with respect to its main diagonal, and $\bs{Z}=\bs{Z}^T$ as well as $\bs{Z}_s=\bs{Z}_s^T$. 
	
	The vectors of currents and voltages associated with this impedance matrix system are partitioned accordingly, i.e., $\bs{i}^T=\begin{bmatrix}{i}_t,\,\bs{i}_s^T,\, i_r\end{bmatrix}$ and $\bs{v}^T=\begin{bmatrix}{v}_t,\,\bs{v}_s^T,\, v_r\end{bmatrix}$ where $\bs{i},\bs{v}\in \mathbb{C}^{N+2}$.
	
\subsection{Objective: Power Transfer Efficiency, Power Gain, Antenna Gain, and Bistatic Radar Cross Section}
%Using the aforementioned quantities, 
%%Using the impedance matrix $\bs{Z}$ according to \eqref{eq:Z}, and the currents $\bs{i}$, 
%the power transfer efficiency (PTE) can be formulated from the ratio of received and transmitted power:
%\begin{equation}
%	\eta=\dfrac{P_r}{P_t}
%	=\dfrac{|i_r|^2R_{\mathrm{rx}}}{\Re\bigl(\bs{i}^H\bs{v}\bigr)}
%	=\dfrac{|i_r|^2R_\mathrm{rx}}{\bs{i}^H(\Re\bs{Z})\bs{i}}\,,%{\operatorname{Re}\bigl(\bs{i}^H\bs{v}\bigr)}\;,
%\end{equation}
%where $R_\mathrm{rx}$ is the receiving resistance. Usually, since impedance-matched transmit and receiving antennas are considered, $R_\mathrm{rx}=R_r=\Re z_r$.

Using the aforementioned quantities, 
%Using the impedance matrix $\bs{Z}$ according to \eqref{eq:Z}, and the currents $\bs{i}$, 
 the \textit{Power Transfer Efficiency} (PTE) $\eta$ of the reduced two-port network can be formulated from the ratio of received and transmitted power according to 
\begin{equation}
	\eta=\dfrac{P_r}{P_t}
	=\dfrac{\frac{1}{2}|i_r|^2R_{\mathrm{rx}}}{\frac{1}{2}\Re\bigl(\bs{i}^H\bs{v}\bigr)}
	=\dfrac{|i_r|^2\Re z_r}{\bs{i}^H(\Re\bs{Z})\bs{i}}\,,%{\operatorname{Re}\bigl(\bs{i}^H\bs{v}\bigr)}\;,
	\label{eq:PTE}
\end{equation}
where due to conjugate-matched transmit and receive antennas, %(i.e., $\Gamma_L=\Gamma_S=0$ due to complex-conjugate loading), 
the receiving resistance is $R_\mathrm{rx}=\Re z_r$.
 %$R_\mathrm{rx}$ is the receiving resistance. Usually, since impedance-matched transmit and receiving antennas are considered, $R_\mathrm{rx}=R_r=\Re z_r$.
%The original framework optimizes the power transfer efficiency (PTE) for wireless power transfer (WPT) systems in the near-field \cite{lang17ant}. 
%However, the PTE $\eta$ is equivalent to the \textit{power gain} $G$ of a two-port, where %$G$~\cite{pozar12} in terms of the two-port S-parameters of the system. %: 
%\begin{equation}
%	G=\dfrac{P_L}{P_\mathrm{in}}=\abs{S_{21}}^2\cdot\dfrac{1}{1-\abs{\Gamma_\mathrm{in}}^2}\cdot \dfrac{1-\abs{\Gamma_L}^2}{\abs{1-S_{22}\Gamma_L}^2}\;.
%\end{equation}
Under these conditions, the \textit{Power Gain} $g$ of the reduced two-port network can be expressed in terms of S-parameters of the two-port network and antenna parameters as follows:
\begin{equation}
	%\left.G\right|_{\Gamma_L=0}
	g=\dfrac{P_r}{P_t}=\dfrac{\abs{S_{21}}^2}{1-\abs{S_{11}}^2}=\dfrac{G_\mathrm{tx}'G_\mathrm{rx}}{L_\mathrm{FS}}=\dfrac{G_\mathrm{tx}G_\mathrm{rx}}{L_\mathrm{FS}}\,\dfrac{\sigma_{\!B}}{4\pi d^2}\,.
	\label{eq:PowerGain}
\end{equation}
$G_\mathrm{tx}$ and $G_\mathrm{rx}$ refer to the gains of the transmitter and receiver antenna, respectively, $\sigma_{\!B}$ is the bistatic radar cross section (BRCS) of the reflector, such that the combined transmitter gain is $G_\text{tx}'=G_\text{tx}\sigma_{\!B}/(4\pi d^2)$, and $L_\mathrm{FS}=(4\pi d/\lambda)^2$ is the free-space loss over the distance $d$ at the free-space wavelength $\lambda$. Hence, optimizing the PTE $\eta$ is equivalent to optimizing the power gain $g$, which in turn is equivalent to optimizing the total transmit antenna gain $G_\mathrm{tx}'$ or the BRCS $\sigma_{\!B}$ of the RIS. %, where $\lambda$ is the free-space wavelength. 

%The same can be expressed for all other receiver antennas in sidelobe and backlobe directions. Moreover, since due to the far-field conditions the coupling between the mainlobe, backlobe and sidelobe receivers as well as the transmitter antenna is very weak, the optimum load impedances %and reactance for the receiver antenna do not need to be obtained via optimization: As the distance $d\rightarrow\infty$, the well-known 
%are the well-known conjugate matches of the respective self-impedances: $z_L^\star =R_L^\star+j X_r\rightarrow \overline{z}_r=z_r'-z_r''$, where $z_r$ denotes the receiver self-impedance. For every receiver, this can be used to obtain the aforementioned condition $\Gamma_L=0$.% (and similarly for all sidelobe and backlobe receeivers).%(i.e. conjugate match, $\Gamma_L=0$).

\subsection{Constraints}

\subsubsection{KVL at the Receiver}
%The optimization problem is set up in terms of the currents $\bs{i}=\begin{bmatrix}\bs{i}_t^T,i_b,i_r\end{bmatrix}{\vphantom{\bs{i}_t}}^{\hspace*{-1pt}T}$. % and the impedance matrix of the system $\bs{Z}$. 
On the receiver end, impedance-matched conditions lead to the KVL  
\begin{equation}
	\begin{bmatrix}z_{tr}\,,\,\bs{z}_{sr}^T\,,\,2\Re z_r
	\end{bmatrix}\bs{i}=0\,.
\end{equation}	
And to obtain a well-defined problem, all phases are expressed in reference to the receiver current by setting $\Im i_r=0$ (via constraint or by reducing the vector of variables accordingly).

The entire optimization problem is first formulated in terms of a (real-valued) vector of unknowns, corresponding to the real and imaginary parts of the currents $\bs{c}^T=\begin{bmatrix}\Re\bs{i}^{\,T},\hspace*{1pt} \Im \bs{i}^{\,T}\end{bmatrix}$ (and similarly for the voltage). KVL and the phase condition can be ensured within affine equality constraints $\bs{A}\bs{c}=\bs{b}$.

\subsubsection{Passivity of RIS-Elements}
%As explained in detail in \cite{lang17}, 
The goal of the optimization procedure is to find the optimally tuned reactances by maximizing the received power while maintaining zero real-power constraints at all ports other than the transmitter and the receiver, i.e., at all RIS antenna ports:
\begin{equation}
	P_n={\textstyle\frac{1}{2}}\bs{c}^T\bs{Q}_n\bs{c}=0\,,\quad \forall n\in[2,N+1]\,.
	\label{eq:Pn}
\end{equation}
Each $\textbf{Q}_n$ corresponds to the $n^\text{th}$ real-valued port impedance matrix. They can  be directly formulated from the impedance matrix $\bs{Z}$ \cite{lang17}, noting two typos corrected in \cite{lang22}. However, since each $\bs{Q}_n$ is an indefinite matrix, \eqref{eq:Pn} amounts to $N$ nonconvex quadratic equality constraints.

Requiring zero real power to be transferred via these ports, means that only reactive power can be exchanged. The corresponding optimal reactive termination for each port $n$ (i.e., each RIS element) is found by 
\begin{equation}
	x_n^\star=-\Im\left({v_n^\star}/{i_n^\star}\right),
\end{equation}where $v_n^\star$ is the optimal unloaded voltage, found from the optimal currents and the impedance matrix via $\bs{v}^\star=\bs{Z}\bs{i}^\star$. Note that the superscript stars denote optimal quantities, not complex conjugates.

\subsubsection{Optional Additional Constraints}
Additional constraints could be added to ensure impedance matching. While that is crucial for optimizing transmit array antennas, it is not required here. Both transmit and receive antennas will be designed to be matched anyway. Furthermore, side- and/or backlobe constraints may be added by including additional receiving antennas in the model in the desired directions/locations. This could be interesting for (S-)RIS applications, but is beyond the scope of this paper.

%
%The difficulty arising when optimizing antennas, particularly passively excited arrays, is that the resulting maximum-gain antenna array solution might be severely mismatched to the desired reference impedance, e.g. $Z_0=50\,\Omega$. 
%
%However, here this is not an issue.

%Multiplying both sides of the right condition in \eqref{eq:matchingConstraintR} by the squared amplitude of its respective input current, $\abs{i_\mathrm{in}}^2=\abs{i_2}^2$, it can be turned into a power constraint:
%\begin{equation}
%	P_\mathrm{in}=P_2\in[P_0/\gamma,P_0\gamma]
%	\quad\Leftrightarrow\quad
%	\gamma P_0\geq P_\mathrm{in}\geq {P_0}/{\gamma}\;.
%	\label{eq:matchingConstraintP}
%\end{equation}
%%\end{subequations}
%Thus, the antenna impedance lying within the required range given by the reference impedance and the mismatch factor is equivalent to the antenna's input power lying in the range given by the power dissipated in the reference (generator) impedance and the mismatch factor.

\subsection{Optimization Problem Formulation}
%Using all these constraints and the cost function, the entire optimization problem can be assembled, leading to 
The following simple Qua\-drat\-ically Constrained Quadratic Program (QCQP) is obtained by combining all these constraints and the cost function:%can be formulated in terms of a real valued vector of the currents $\textbf{c}=\left[\bm{i}_t'^{T}, \bm{i}_s', i_r',\bm{i}_t''^{T}, \bm{i}_s'', i_r''\right]$ as a Quadratically Constraint Quadratic Program (QCQP) of the form
\begin{equation}
		\begin{aligned}
		%p^\star:\quad
		P_{t,\mathrm{min}}^\mathrm{QCQP}\!=\min_{\bs{c}}\;&P_\mathrm{in}={\textstyle\frac{1}{2}}\bs{c}^T({\textstyle\sum_n\!\bs{Q}_n})\bs{c}\\
		\mathrm{s.t.}\;
		&P_n={\textstyle\frac{1}{2}}\bs{c}^T\bs{Q}_n\bs{c}=0\,,
		&& n\in[2,N+1]\\[0.25em]
		&\bs{A}\bs{c}=\bs{b}\;.
		%[z_{tr},\bs{z}_{sr}^T,2\Re z_r]\bs{c}=\bs{b}\;.
		%&c_N=\sqrt{2/R_L}
		%&P_\mathrm{rx}={\textstyle\frac{1}{2}}\bs{c}^T\bs{R}\bs{c} = 1\;.
	\end{aligned}
	\label{eq:qcqp2}
\end{equation}

%Furthermore, only dipole $n=2$ is excited actively; the dipoles $n\neq2$ (including the receiver) are only excited indirectly via mutual coupling. Thus, there is neither insertion nor extraction of real power via those ports, i.e.
%\begin{equation}
%	P_n=\dfrac{1}{2}\bs{c}^T\bs{Q}_n\bs{c}=0\quad\text{for }n\neq2\;,
%	\label{eq:passivity_constraint_c}
%\end{equation}
%where $\bs{Q}_n=\bs{Q}_n^H$ is the real-valued form of the loaded port impedance matrix (i.e. the part of the loaded impedance matrix corresponding to real power flow through the $n$th port) \cite{lang17b}. Note the correction in the appendix. 

However, this QCQP %program \eqref{eq:qcqp2} 
is \textit{nonconvex}. It belongs to the class of problems known to be notoriously difficult to solve. 
Such problems are frequently referred to as ``unsolvable'', because the computational cost can quickly exceed any feasible means, at least for a large number of variables.	

\subsection{Semidefinite Relaxation (SDR)}
From the nonconvex QCQP \eqref{eq:qcqp2}, a Semidefinite Program (SDP) \cite{boyd04} can be formulated in terms of the current matrix $\bs{C}=\bs{c}\bs{c}^T$ using the cyclic permutation property of the trace $\bs{c}^T\bs{Q}\bs{c}=\tr(\bs{c}^T\bs{Q}\bs{c})=\tr(\bs{Q}\bs{c}\bs{c}^T)=\tr(\bs{Q}\bs{C})$, resulting in:
\begin{equation}
	\begin{aligned}
		P_{t,\mathrm{min}}^\mathrm{SDR}=\min_{\bs{c},\bs{C}}\;&{\textstyle\frac{1}{2}}\tr({\textstyle\sum_n\!\bs{Q}_n}\bs{C})\\[-0.15em]
		\mathrm{s.t.}\;
		&\tr(\bs{Q}_n\bs{C}) = 0\,, && n\in[2,N+1]\\[0.25em]
		&\bs{A}\bs{c} = \bs{b}\\
		&\begin{bmatrix}\bs{C}\!&\!\bs{c}\\\bs{c}^T\!&\!1\end{bmatrix}\succeq0\;.
	\end{aligned}
	\label{eq:sdp}
\end{equation}
%Note that the fact that \eqref{eq:sdp} is an SDP is well visible in the last, typical, semidefinite constraint.

The programs \eqref{eq:qcqp2} and \eqref{eq:sdp} are in general not completely equivalent. 
That would require that $\bs{C}-\bs{c}\bs{c}^T=0$. 
However, this condition is equivalent to {both} $\bs{C}-\bs{c}\bs{c}^T\succeq0$ \textit{and} $\bs{C}-\bs{c}\bs{c}^T\preceq0$. 
The elimination of the second of these two constraints is the typical semidefinite relaxation applied to obtain a solvable convex SDP. Such programs are reliably and efficiently solvable, even for a (reasonably) large number of variables, using dedicated algorithms, such as SDPT3 \cite{tutuncu03}.

%The last constraint is typical for such SDPs. %The SDP \eqref{eq:sdp} can be solved reliably and efficiently using dedicated algorithms, such as SDPT3 \cite{tutuncu03}. 
% In Matlab, it can be implemented comfortably using CVX \cite{cvx,grant08} and solved for example using the standard solver SDPT3 \cite{toh99,tutuncu03}. %To enhance numerical precision, redundant equality constraints may be added, as discussed in \cite{lang16}.

\subsection{Tightness of the Problem}
A test has to confirm whether the optimum solution of the SDP \eqref{eq:sdp} actually coincides with the optimum of the original QCQP \eqref{eq:qcqp2}. %; in such cases, the relaxation is referred to as \textit{tight}. 
To check for this desired property, the following (relative) \textit{tightness error} (TE) \cite{lang17}, is used:
\begin{equation}\label{eq:tightness_error}
	\epsilon={\norm{\bs{C}^*-\bs{cc}^T}}/\left({\bs{c}^T\bs{c}}\right)\,.
\end{equation}
If $\epsilon$ is sufficiently small, the relaxation was \textit{tight} and the solution of the SDP \eqref{eq:sdp} corresponds to the solution of the original nonconvex QCQP \eqref{eq:qcqp2}. 

However, if that is not the case, the obtained solution is infeasible with respect to the original problem, because $\bs{C}$ cannot be factored into $\bs{c}\bs{c}^T$. %, i.e., the solution cannot be obtained by a single set of currents $\bs{c}$. 
It may still serve as a good starting point for another algorithm to find a suitable "next-best" feasible solution. However, it cannot be proven that the solution found in this way is actually the optimal solution to the original problem.%nonconvex QCQP \eqref{eq:qcqp2}.

\vspace*{4pt}
\section{Simulation Examples}
\begin{figure*}[!htb]\centering
		\hspace*{-0.2mm}\includegraphics{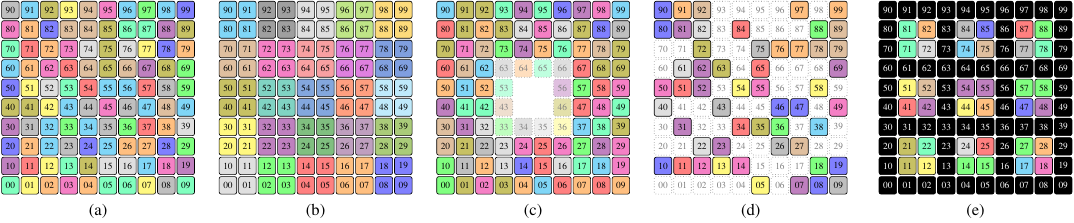}
		\vspace*{-6.5mm}
	\caption{Configurations of Sparse RIS (SRIS) considered in this paper: (a) fully tunable RIS, (b) parallel-connected $2\!\times\!2$ clusters, (c) $2\!\times\!2$ or $4\!\times\!4$ center removed, (d) a random selection of 50\% or 75\% of the elements, and (e) nine $2\!\times\!2$-subregions with the other antennas short-circuited.}
	\label{fig:configurations}
\end{figure*}
\subsection{Setup}
The simple $\text{10}\!\times\!{10}$-element printed dipole surface, shown in \figref{fig:hfss_setup}, serves as illustrative test case. The dipoles are modeled using finite-conductivity sheets (copper conductivity) to account for losses. They are uniformly separated by $\lambda/2$ at the operating frequency of 3.75\,GHz and lie on a 1.6\,mm thick FR4 substrate backed by a layer of perfect conductor. 
%The total dimensions of the SRIS is 40\,cm by 40\,cm. 

The transmit and receive antennas are two standard pyramidal horn antennas with about 16\,dB gain. %This surface may not generally look like many of the current designs, however, that really doesn't matter.
 Their distance to the RIS is always $d=2D^2/\lambda$ (with $D$ being the diagonal of the surface), such that far-field conditions apply. However, this is not a requirement for the method.

Note: While it is certainly computationally expensive to obtain the $102\times102$ impedance matrix of the entire system, such simulations are often carried out at some point during the design of such a surface. Moreover, approximative methods could be used in the far-field.

\subsection{SRIS Configurations}
The following four SRIS configuration types, depicted in \figref{fig:configurations}(b)--(e), are compared to the full RIS, shown in \figref{fig:configurations}(a):
\begin{itemize}
	\item parallel-connected $2\!\times\!2$-element clusters,
	\item removing a $2\!\times\!2$- or a $4\!\times\!4$-element center,
	\item randomly enabling 75\% or 50\% of all elements,
	\item nine $2\!\times\!2$-element subarrays, the other elements shorted.
\end{itemize}
The case with all RIS elements open-connected, i.e., predominantly a perfectly conducting square reflector with an RCS of about $17$\,dB in the normal direction, serves as base reference.

%\balance

%For the limited scope of this paper, only the angle $\alpha$ is swept in this study.

%\section{Simulation Results}
%\subsection{Setup}

\subsection{Results}

\subsubsection{Tightness} 
\begin{figure}[!htb]\centering
	\vspace*{-1mm}
	\begin{tikzpicture}\small\sffamily
		\node at(0,0){\includegraphics[width=0.95\columnwidth,clip,trim={6mm 0 6mm 1mm}]{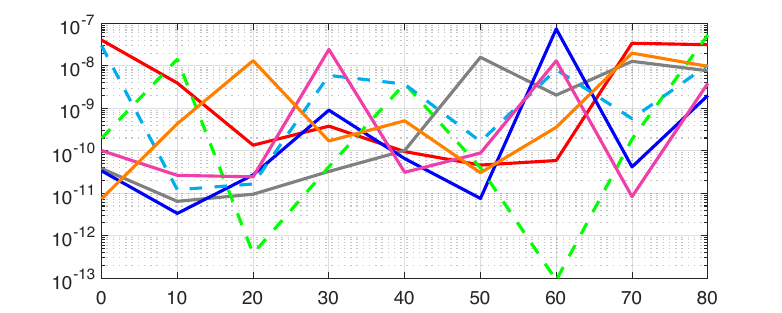}};
		\node[rotate=90,scale=0.85] at(-4.35,0){Tightness error $\epsilon$};
		\node[scale=0.85] at(0.25,-2.1){Angle off broadside (degrees)};
	\end{tikzpicture}
	\captionsetup{justification=raggedleft}
	\vspace*{-6mm}
	\caption{These small tightness errors confirm that the relaxations are tight in all cases and that the optimum of the nonconvex QCQP \eqref{eq:qcqp2} has been found at every point. See \figref{fig:results} for a legend.}
	\label{fig:tightness}
	\vspace*{-2mm}
\end{figure}

First of all, the graphs in \figref{fig:tightness} reveal that the relaxation was tight in all cases. Thus, indeed, the optimum solution of the QCQP~\eqref{eq:qcqp2} has been found at every point. %its optimum solution coincides with the optimum solution of the original nonconvex problem.

\subsubsection{Optimized BRCS} 
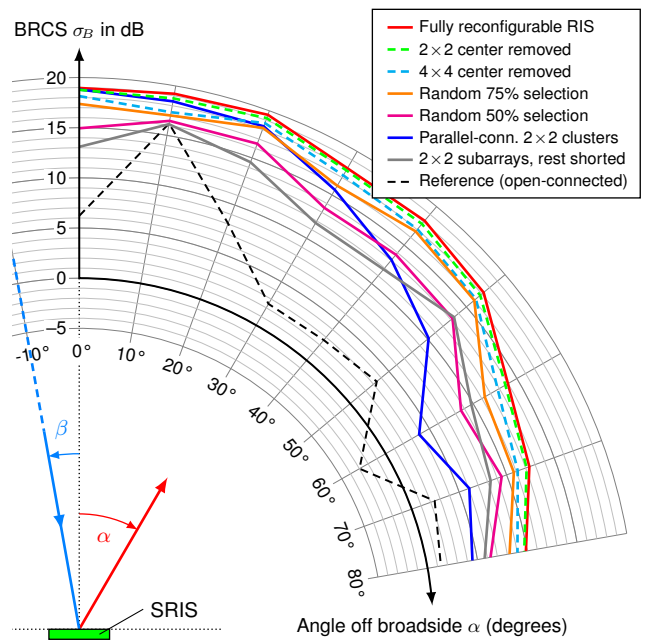
\begin{figure}[!b]\centering
	\vspace*{-4mm}
	\begin{tikzpicture}[line width=0.5pt,>=latex,scale=1.33]\sffamily\small
		
		\clip (-0.66,-0.12)rectangle(5.60,6.2);
		\def\rad{3.5}
		
		\draw[densely dotted] (0,0)--(0,\rad);
		\draw[densely dotted] (0.4*\rad,0)--(-0.4*\rad,0);

		\foreach \al in {-10,10,20,...,80}{
			\draw[gray,line width=0.35pt] (90-\al:\rad-0.6)--(90-\al:\rad+2);
			\node[rotate=-\al,scale=0.75] at(90-\al:\rad-0.3-0.42){\al{}°\!};
		}
		\def\al{0}
			\node[rotate=-\al,scale=0.75,fill=white] at(90-\al:\rad-0.3-0.42){\al{}°\!};
		\foreach \dB in {-0.4,-0.3,...,1.9}{
			\draw[lightgray,line width=0.25pt] (110:\rad+\dB)arc(110:10:\rad+\dB);
		}
		\foreach \dB in {-0.5,0,...,2}{
			\draw[gray,line width=0.35pt] (110:\rad+\dB)arc(110:10:\rad+\dB);
		}
		%\draw[line width=1pt] (0,\rad)arc(90:0:\rad);
		\node[anchor=east,scale=0.75,fill=white,inner sep=1pt] at(-0.1,\rad-0.5){--5};
		\node[anchor=east,scale=0.75,fill=white,inner sep=1pt] at(-0.1,\rad){0};
		\node[anchor=east,scale=0.75,fill=white,inner sep=1pt] at(-0.1,\rad+0.5){5};
		\node[anchor=east,scale=0.75,fill=white,inner sep=1pt] at(-0.1,\rad+1.0){10};
		\node[anchor=east,scale=0.75,fill=white,inner sep=1pt] at(-0.1,\rad+1.5){15};
		\node[anchor=east,scale=0.75,fill=white,inner sep=1pt] at(-0.1,\rad+2.0){20};
		
		\draw[red,->,line width=1pt] (0,0)--(60:\rad*0.5);
		\draw[blue!50!cyan,line width=1pt,densely dashed] (90+10:\rad+0.25)--(90+10:\rad-1.5);
		\draw[blue!50!cyan,line width=1pt] (90+10:\rad-1.5)--coordinate(temp)(0,0);
		\draw[blue!50!cyan,line width=1pt,<-] (temp)--+(90+10:2pt);
		\draw[fill=green] (0,0)+(-0.3,0)rectangle(0.3,-0.1);
		\draw[thin] (0.2,-0.05)--+(30:0.5)node[right,scale=0.85]{SRIS};
		%\node[scale=0.9] at(0,-0.3){SRIS};
		\draw[red,->] (0,1.15)arc(90:63:1.15)--+(60-90:2pt);
		\node[red,scale=0.9] at(75:0.95){$\alpha$};
		\draw[blue!50!cyan,->] (0,1.75)arc(90:98:1.75)--+(97+90:2pt);
		\node[blue!50!cyan,scale=0.9] at(95:2){$\beta$};
		
		% full RIS
		\draw[red,line width=1pt] 
		 (90-0:\rad+18.9623/10)--(90-10:\rad+19.2085/10)--(90-20:\rad+19.5764/10)--(90-30:\rad+18.0635/10)--(90-40:\rad+18.1643/10)--(90-50:\rad+17.2252/10)--(90-60:\rad+13.9703/10)--(90-70:\rad+12.3876/10)--(90-80:\rad+9.8708/10);
		
		% reference (open-connected)
		\draw[black,densely dashed,line width=0.75pt] 
		 (90-0:\rad+6.2150/10)--(90-10:\rad+16.1921/10)--(90-20:\rad+8.0503/10)--(90-30:\rad+2.4495/10)--(90-40:\rad+2.6329/10)--(90-50:\rad+3.4223/10)--(90-60:\rad-2.9822/10)--(90-70:\rad+2.4205/10)--(90-80:\rad+1.3330/10);
	
	% 2x2 subarrays shorted
	\draw[blue,line width=1pt] 
		 (90-0:\rad+18.8128/10)--(90-10:\rad+18.4367/10)--(90-20:\rad+18.3567/10)--(90-30:\rad+15.5762/10)--(90-40:\rad+13.0632/10)--(90-50:\rad+10.1210/10)--(90-60:\rad+3.8318/10)--(90-70:\rad+6.0671/10)--(90-80:\rad+4.4983/10);
	
		% no center 1
	\draw[green,densely dashed,line width=1pt] 
		(90-0:\rad+18.7756/10)--(90-10:\rad+18.7372/10)--(90-20:\rad+19.23/10)--(90-30:\rad+17.8374/10)--(90-40:\rad+17.7513/10)--(90-50:\rad+16.7878/10)--(90-60:\rad+13.6823/10)--(90-70:\rad+12.0962/10)--(90-80:\rad+9.7057/10);
		% (90-0:\rad+18.8630/10)--(90-10:\rad+18.9515/10)--(90-20:\rad+19.2857/10)--(90-30:\rad+17.9866/10)--(90-40:\rad+17.9599/10)--(90-50:\rad+16.9997/10)--(90-60:\rad+13.7981/10)--(90-70:\rad+12.1896/10)--(90-80:\rad+9.8536/10); 
		
		% no center 2
	\draw[cyan,densely dashed,line width=1pt] 
		(90-0:\rad+18.1511/10)--(90-10:\rad+17.3363/10)--
		(90-20:\rad+18.6977/10)--(90-30:\rad+17.2442/10)--
		(90-40:\rad+17.0578/10)--(90-50:\rad+16.2266/10)--
		(90-60:\rad+12.8087/10)--(90-70:\rad+11.1655/10)--
		(90-80:\rad+9.0232/10);
		
   % random 75%
	\draw[orange,line width=1pt] 
		(90-0:\rad+17.3707/10)--(90-10:\rad+17.0091/10)--
		(90-20:\rad+18.1852/10)--(90-30:\rad+16.0561/10)--
		(90-40:\rad+16.7781/10)--(90-50:\rad+16.0194/10)--
		(90-60:\rad+11.3269/10)--(90-70:\rad+10.7528/10)--
		(90-80:\rad+8.2302/10);

	% random 50%
	\draw[magenta,line width=1pt] 
		(90-0:\rad+14.9440/10)--(90-10:\rad+16.4602/10)--
		(90-20:\rad+16.5120/10)--(90-30:\rad+13.5059/10)--
		(90-40:\rad+13.7437/10)--(90-50:\rad+13.2154/10)--
		(90-60:\rad+8.6173/10)--(90-70:\rad+9.4697/10)--
		(90-80:\rad+6.3159/10);

	\draw[gray,line width=1pt] 
		 (90-0:\rad+13.0939/10)--(90-10:\rad+16.1260/10)--(90-20:\rad+14.5471/10)--(90-30:\rad+11.7680/10)--(90-40:\rad+11.8376/10)--(90-50:\rad+13.4107/10)--(90-60:\rad+9.7889/10)--(90-70:\rad+8.3120/10)--(90-80:\rad+5.7276/10);
		
		\draw[line width=0.75pt,<->] (0,\rad+2.3)node[above=1pt,scale=0.8]{BRCS $\sigma_{\!B}$ in dB}--(0,\rad)arc(90:5:\rad)--+(5-90:2pt)node[below=2pt,scale=0.8]{Angle off broadside $\alpha$ (degrees)};
		
		\begin{scope}[xshift=3.05cm,yshift=6.0cm,every node/.append style={scale=0.77}]\footnotesize
			\def\d{0.22}
			\draw[fill=white] (-0.15,0.18)rectangle(2.5,-1.7);
			\draw[red,line width=1pt] (0,0)--+(0.25,0)node[right,text=black]{Fully reconfigurable RIS};
			\draw[green,densely dashed,line width=1pt] (0,-\d)--+(0.25,0)node[right,text=black]{$\text{2}\!\times\!\text{2}$ center removed};
			\draw[cyan,densely dashed,line width=1pt] (0,-2*\d)--+(0.25,0)node[right,text=black]{$\text{4}\!\times\!\text{4}$ center removed};
			\draw[orange,line width=1pt] (0,-3*\d)--+(0.25,0)node[right,text=black]{Random 75\% selection};
			\draw[magenta,line width=1pt] (0,-4*\d)--+(0.25,0)node[right,text=black]{Random 50\% selection};
			\draw[blue,line width=1pt] (0,-5*\d)--+(0.25,0)node[right,text=black]{Parallel-conn. $\text{2}\!\times\!\text{2}$ clusters};
			\draw[gray,line width=1pt] (0,-6*\d)--+(0.25,0)node[right,text=black]{$\text{2}\!\times\!\text{2}$ subarrays, rest shorted};
			\draw[black,densely dashed,line width=0.75pt] (0,-7*\d)--+(0.25,0)node[right,text=black]{Reference (open-connected)};
		\end{scope}
		
	\end{tikzpicture}
	\vspace*{-1mm}
	\caption{The resulting bistatic radar cross section (BRCS) of the RIS, calculated from the maximum PTE values, for different SRIS configurations and angles $\alpha=0\,...\,80$° off broadside; constant angle of incidence $\beta=-10$°.}
	\label{fig:results}
	\vspace*{-1mm}
\end{figure}

Using the equivalence between \eqref{eq:PTE} and \eqref{eq:PowerGain}, the bistatic radar cross section is calculated from the obtained maximum PTE values. They are displayed vs. the angle off broadside to the RIS $\alpha$ in \figref{fig:results},%, 
when the angle of incidence (transmit antenna) is constant $\beta=-10$°.

As can be seen by the dashed green line, removing the four center elements has negligible effect, the larger center (dashed blue) reduces the BRCS by 1-2\,dB. Saving these additional 10\% of reconfigurable elements comes at a (small) price.

Reducing the SRIS randomly by 25\% of its reconfigurable elements has a larger effect, but 
%The SRIS resulting by reducing the RIS randomly by 25\% of its elements
the SRIS can still provide a significantly larger BRCS, also at angles far off broadside. Even with only 50\% reconfigurable elements, the SRIS still manages to provide a substantially improved BRCS compared to the reference reflector, which only provides the specular reflection at +10°, as expected.

The SRIS that results from parallel-connected $2\!\times\!2$ clusters, i.e., essentially an RIS of similar size but with only 25 instead of 100 reconfigurable elements, performs very well for small angles of reflections and only at angles beyond 45° starts falling behind. Lastly, 36 reconfigurable elements in nine $2\!\times\!2$-subarrays can provide better performance at flat angles, but falls behind at steep angles.

There is no clear winner at this point, but a few of these strategies definitely deserve further investigation.

\vspace*{4pt}
\section{Summary \& Conclusion}
First, this paper introduces the concept of Sparse Reconfigurable Intelligent Surfaces (SRIS): While RIS could play a crucial role in future wireless networks to enhance both coverage and energy efficiency, a large number of reconfigurable elements (e.g., varactor diodes as tunable capacitors) may prevent the second reason. Reducing the number of reconfigurable elements clearly affects the overall performance, while improving energy efficiency. This trade-off needs to be studied carefully, to find an optimal consensus.

Second, the versatile maximum power transfer optimization framework is adapted and applied to SRIS configurations to find their optimal performance and provide helpful insight for finding the aforementioned  consensus between number and placement of reconfigurable elements and the performance of the SRIS. As shown, the involved relaxation is tight in all these cases, which means that the global maximum in terms of performance, measured via bistatic radar cross section (BRCS), and the corresponding optimal loading reactances could be found each time.

Finally, preliminary results show that different configurations to achieve the desired sparsity have different effects on reducing the number of reconfigurable components and the resulting RIS performance. This opens the way for further investigations, possibly also including constraints on maximum reflection levels in undesired directions.

%\subsection{Balance Last Page}
%To balance the columns on the last page, the package ``balance.sty'' can be used. Alternativey, ``IEEEtran.cls'' provides the commands ``\textbackslash IEEEtriggeratref'' and ``\textbackslash IEEEtriggercmd'' (see the documentation in this \LaTeX\ template).


\section*{Acknowledgment}

This work was carried out as part of the collaborative CHIST-ERA project \textit{"Towards Sustainable ICT: Sparse Ubiquitous Networks based on Reconfigurable Intelligent SurfacEs (SUNRISE)"}. The authors of OST are supported by the Swiss National Science Foundation (SNF) with Grant 203784. %The author from UCD is supported by...

% \IEEEtriggeratref{7}
% \IEEEtriggercmd{\enlargethispage{-20cm}}
%\bibliographystyle{IEEEtran}
%Unless there are six authors or more give all authors’ names; do not use ``\textsl{et al.}''. Papers that have not been published, even if they have been submitted for publication, should be cited as ``unpublished'' \cite{Elissa}. Papers that have been accepted for publication should be cited as ``in press'' \cite{Nicole}. Capitalize only the first word in a paper title, except for proper nouns and element symbols.

\bibliographystyle{myIEEEtran}

\end{document}